\documentclass[11pt, sort]{elsarticle}
\usepackage{graphicx}
\usepackage{subcaption}
\usepackage{natbib}
\usepackage{hyperref}
\usepackage{mathtools}
\usepackage{amsfonts}
\usepackage{siunitx}
\usepackage{booktabs}
\usepackage{multirow}
\usepackage{titlesec}
\usepackage{caption}
\usepackage{xcolor}
\usepackage{tikz}
\usetikzlibrary{arrows.meta}
\usepackage[margin=3cm]{geometry}

\mathtoolsset{showonlyrefs}

\hypersetup{
    colorlinks = true,
    urlcolor   = blue,
    citecolor  = black,
}

\newcommand{\RomanNumeralCaps}[1]

\titleformat*{\section}{\normalfont\fontsize{12}{15}\bfseries}
\titleformat*{\subsection}{\normalfont\fontsize{11}{15}\bfseries}

\DeclareMathOperator{\tr}{tr}
\DeclareMathOperator{\dev}{dev}

\newcommand{\bA}{\boldsymbol{\mathrm{A}}}

\newcommand{\bD}{\boldsymbol{\mathrm{D}}}
\newcommand{\bI}{\boldsymbol{\mathrm{I}}}
\newcommand{\bS}{\boldsymbol{\mathrm{S}}}
\newcommand{\bO}{\boldsymbol{\mathrm{O}}}

\newcommand{\bsigma}{\boldsymbol{\sigma}}

\newcommand{\btheta}{\boldsymbol{\theta}}

\newcommand{\dd}{\mathrm{d}}
\newcommand{\bu}{\boldsymbol{u}}
\newcommand{\bv}{\boldsymbol{v}}
\newcommand{\bzero}{\boldsymbol{0}}

\newcommand{\RR}{\mathbb{R}}
\newcommand{\Sbb}{\mathbb{S}}
\newcommand{\Obb}{\mathbb{O}}

\newcommand{\bff}{\boldsymbol{f}}

\newcommand{\bn}{\boldsymbol{n}}

\newcommand{\Cc}{\mathcal{C}}

\newcommand{\Jc}{\mathcal{J}}
\newcommand{\Sc}{\mathcal{S}}
\newcommand{\Vc}{\mathcal{V}}

\begin{document}
\begin{frontmatter}
\title{Learning parameter-dependent shear viscosity from data, \\ with application to sea and land ice}
\author{Gonzalo G.~de Diego}
\author{Georg Stadler}
\affiliation{organization={Courant Institute of Mathematical Sciences, New York University},
            country={USA}
            }

\begin{abstract}
    Complex physical systems which exhibit fluid-like behavior are often modeled as non-Newtonian fluids. A crucial element of a non-Newtonian model is the rheology, which relates inner stresses with strain-rates. We propose a framework for inferring rheological models from data that represents the fluid's effective viscosity with a neural network. By writing the rheological law in terms of tensor invariants and tailoring the network's properties, the inferred model satisfies key physical and mathematical properties, such as isotropic frame-indifference and existence of a convex potential of dissipation. Within this framework, we propose two approaches to learning a fluid's rheology: 1) a standard regression that fits the rheological model to stress data and 2) a PDE-constrained optimization method that infers rheological models from velocity data. For the latter approach, we combine finite element and machine learning libraries. We demonstrate the accuracy and robustness of our method on land and sea ice rheologies which also depend on external parameters. For land ice, we infer the temperature-dependent Glen's law and, for sea ice, the concentration-dependent shear component of the viscous-plastic model. For these two models, we explore the effects of large data errors. Finally, we infer an unknown concentration-dependent model that reproduces Lagrangian ice floe simulation data. Our method discovers a rheology that generalizes well outside of the training dataset and exhibits both shear-thickening and thinning behaviors depending on the concentrations. 
\end{abstract}


\end{frontmatter}

\section{Introduction}

Non-Newtonian fluid models are often used as continuum representations of complex physical systems which, over suitable spatial and temporal scales, display a fluid-like behavior. Examples of these complex systems include polymer solutions \cite{larson2013}, colloidal suspensions \cite{Guazzelli18}, granular flows \cite{jop2006}, and Earth's mantle \cite{ranalli1995}. Two additional examples explored in this article are ice sheets (i.e.~land ice) \cite{greve2009} and sea ice \cite{feltham2008}. 

Non-Newtonian models are generally written as a system of partial differential equations (PDEs) that consist of conservation laws and constitutive models. The constitutive model for a fluid, which is also referred to as its rheology, relates inner stresses with the deformation of the material. Traditionally, rheological models have been built in two ways: by finding a function with a simple analytical expression that fits the data (phenomenological models) or by deriving the model from first principles using e.g.~kinetic theory \cite{jenkins1983, toppaladoddi2025}. 
Finding accurate models via these two approaches is often a difficult task due to the inherent complexity of certain materials. 

Data-driven methods are an alternative approach to finding accurate constitutive models for complex materials.
These methods integrate data with mathematical models to infer physically-consistent material laws.
One such class of methods assumes a predefined functional form for the rheology that contains unknown scalar parameters. Typically, gradient-based optimization is used to choose these parameters in order to fit the observation data. In doing so, the governing equations are either enforced as constraints \cite{hu2024,akerson2025, kontogiannis2025}, or their residuals are penalized as done in physics-informed neural networks (PiNNs), \cite{wang2025, lardy2025}.

Using a predefined functional form for the rheology requires knowledge about the material, which may not be available. Recently, several methods have been proposed to avoid using predefined functional forms through the application of scientific machine learning (sciML) techniques. In \cite{saadat2022, mahmoud2024, thakur2024}, for instance,
rheologies are identified within a library of functions using sparsification techniques that promote simple models. 
Alternatively, and of particular importance to this work, a class of methods represent the rheology with neural networks (NNs) \cite{reyes2021, bruer2022, lennon2023, parolini2025, sunol2025, dediego2025}. 
The availability of powerful machine learning libraries with automatic differentiation (AD) capabilities and the approximation properties of NNs make this approach particularly appealing. Before discussing details of this related work, we summarize our approach.

\subsection{Approach and related literature}
In this work, we present a framework for data-driven inference of rheological models from stress and velocity data for non-Newtonian viscous fluids. While we represent rheological functions with an NN, we rely on well established (finite element) discretizations for the governing differential equations. By characterizing the rheological model in terms of tensor invariants to guarantee isotropic frame-indifference, NNs represent unambiguous physical quantities, yielding a fully interpretable framework. For example, in the one-dimensional and incompressible setting explored in this work, the NN represents the effective shear viscosity. 

Next, we review other contributions that also represent rheological models with NNs. They differ in e.g.~what kind of data the NNs are trained on (stress or velocity data) or what physical and mathematical principles are enforced in the model.   
In \cite{lennon2023},  viscoelastic rheologies are represented with NNs and trained with experimental polymer data. The NNs are trained on stress data by solving the time-dependent initial value problem governing a viscoelastic rheology. In \cite{parolini2025}, convex NNs are fitted to stress data via standard regression to represent the rheology for a set of shear-thinning fluids. Convexity of the NN is enforced to ensure well-posedness of the underlying PDEs. A class of monotone rheologies built with an NN are fitted to both rheology and velocity data in \cite{sunol2025}. This approach is applied to synthetic data to recover known rheological models. Finally, in \cite{dediego2025}, the authors of this article apply the method described here to infer an unknown rheology that reproduces synthetic sea ice Lagrangian floe data; we summarize these results briefly in section \ref{subsec:DEM} below. 
This article represents an extension of \cite{dediego2025} in that we derive key physical and mathematical principles that the rheology should satisfy in two-dimensions, we present a streamlined implementation applicable to a wider class of fluids, and we explore the method's accuracy and sensitivity to data errors in synthetic settings in land and sea ice models for which the rheology is known.

\subsection{Contributions and limitations}
The main contributions of this work are as follows:
(1) We derive and discuss the incorporation of three critical physical and mathematical principles for the rheological model: isotropy, the existence of an underlying potential of dissipation, and the second law of thermodynamics. The rheological models satisfying these principles are more general than the class of rheologies considered in \cite{parolini2025, sunol2025}. 
(2) We treat the governing equation during training using methods from PDE-constrained optimization, which are combined with NN libraries.
This should be contrasted with e.g.~PiNNs-based methods, as in \cite{thakur2024, wang2025}, where the sum of the data misfit and the PDE residual is minimized, which may affect the numerical PDE solution accuracy in the presence of large data errors. 
(3) We illustrate numerically that our framework performs robustly for velocity and stress data that include large errors. In particular, we infer two known rheological models (Glen's law for land ice, and the effective shear viscosity of the viscous-plastic model for sea ice) from noisy data. To demonstrate our ability to discover unknown models, we infer a rheology from sea ice data generated with a discrete element method (DEM) that tracks individual ice floes.

The limitations of this work are:
(1) We confine ourselves to the inference of a parameter-dependent effective shear viscosity function. This function fully defines the constitutive relation for incompressible models, as used for land ice, while it only partially characterizes this relationship for compressible flows, as utilized in the context of sea ice.
(2) The inferred rheology functions involve only a single external parameter (concentration for sea ice, and temperature for land ice). However, our approach easily extends to multiple external parameters. 
(3) Our implementation relies on an external coupling between a finite element library for solving PDEs and a machine learning library for parameterizing NN functions. This coupling is not straightforward, particularly when it comes to composing the differentiation capabilities of both libraries. One solution to this is to develop PDE solvers in an AD software, as in \cite{sunol2025, alhashim2025}, and to rely on AD capabilities also for the governing PDE rather than on adjoints.

\section{A framework for rheology inference}\label{sec:framework}
Let $\Omega \subset \RR^2$ be the domain occupied by a viscous fluid. Denoting the fluid velocity by $\bu$, the Cauchy stress tensor by $\bsigma$, and the (possibly nonlinear) external forces by $\bff$, momentum balance is given by
\begin{align}\label{eq:momentum_balance}
    \rho\frac{\bD\bu}{\bD t} -  \nabla\cdot \bsigma = \bff(\bu) \quad \text{in $\Omega$}.
\end{align}
Here, $\rho$ is the density of the fluid and $\bD\bu/\bD t$ the material time derivative of the velocity. To build a closed system of equations, a rheological model
\begin{align}\label{eq:rheology}
    \bsigma = \bsigma(\bD\bu)
\end{align}
that relates the Cauchy stress tensor $\bsigma$ to the strain-rate tensor 
\begin{align}
    \bD\bu := \frac{1}{2}\left( \nabla\bu + \nabla\bu^\top\right)
\end{align}
is required. Here, we present a framework for inferring  function \eqref{eq:rheology} from data. Below, in section \ref{subsec:isotropy}, we constrain the functional form of \eqref{eq:rheology} to ensure that the resulting model is isotropic, determined by a convex potential and dissipative. Then, in section \ref{subsec:simplified}, we restrict our attention to simplified configurations that emerge whenever the problem is one-dimensional or the fluid is incompressible; these are the  configurations that will be considered for the remainder of the article. Then, in sections \ref{subsec:shear_visc_as_NN} and \ref{subsec:optimizing_the_NN}, we write the rheological model in terms of an NN and introduce two loss functions (with stress data or with velocity data) whose optimization results in the NN weights. Finally, section \ref{subsec:comp_setup} details the computational implementation. 

\subsection{Isotropic rheology in two dimensions}\label{subsec:isotropy}

In formulating a general functional form for the rheology of a viscous fluid, we take as a starting point an explicit and local dependence between the Cauchy stress tensor $\bsigma$ and the strain-rate tensor $\bD\bu$:
\begin{align}\label{eq:general_rheology}
	\bsigma = -p\bI + \Cc(\bD\bu)
\end{align}
for a function $\Cc:\Sbb_2\to\Sbb_2$, where $\Sbb^2$ denotes the set of symmetric matrices in $\RR^{2\times 2}$. The assumptions regarding explicitness and locality are a constraint that simplifies both our theoretical and computational framework. Relaxing this constraint could be an intriguing subject for future research.
In the context of our application to sea ice, locality may be the most restricting of the two assumptions, since granular flows are known to exhibit nonlocal behavior \cite{kamrin2012}. For further information on implicit rheologies, we refer the reader to \cite{rajagopal2006}. 

We now proceed to characterize and restrict the functional form of $\Cc$. We first note that the pressure $p$ in \eqref{eq:general_rheology} represents an equilibrium pressure, such that
\begin{align}\label{eq:defn_eq_pressure}
    \bsigma = -p\bI \quad \text{when $\bD\bu = 0$}.
\end{align}
Therefore, we must have that $\Cc(\bD\bu) = 0$ when $\bD\bu = 0$. The pressure $p$ and the function $\Cc$ depend on additional scalar fields that, in fluid dynamics, represent thermodynamic quantities such as density and temperature. In the applications explored in sections \ref{sec:land_ice} and \ref{sec:sea_ice}, we consider the dependence on an additional scalar field that represents the temperature, when applied to land ice, or the sea ice concentration, when the fluid represents sea ice. However, we ignore dependence throughout this section to simplify notation. 

Under the assumption that $\Cc$ is \textit{isotropic}, this function can be characterized more precisely. Given the set of orthogonal matrices,
%
	$\Obb_2 := \left\lbrace \bO\in\RR^{2\times 2} : \bO\bO^\top = \bO^\top\bO = \bI \right\rbrace$
%
we say that the function $\Cc$ is isotropic if
\begin{align}
	\Cc(\bO\bA\bO^\top) = \bO \Cc(\bA)\bO^\top
\end{align}
holds for all $\bO\in\Obb_2$. It is possible to show (see \cite[Section 37]{gurtin1982}) that $\Cc$ is isotropic if and only if there exist two scalar functions $\psi_i:\RR\times \RR_{\geq 0}\to\RR$ for $i = 1, 2$ such that 
\begin{align}\label{eq:representation_of_C}
	\Cc(\bA) = \psi_1(\iota_{\bA}) \bI + \psi_2(\iota_{\bA}) \dev{\bA},
\end{align}
where $\iota_{\bA}\in\RR \times \RR_{\geq 0}$ denotes the principal invariants of the matrix $\bA\in\Sbb_2$ and $\dev{\bA}$ the deviatoric part of $\bA$, given by
\begin{align}
    \dev{\bA} := \bA - \frac{1}{2}(\tr{\bA})\bI,
\end{align}
with $\tr{\bA}$ denoting the trace of $\bA$. For convenience, we represent the principal invariants $\iota_{\bA}$ of a symmetric matrix $\bA\in\Sbb_2$ as 
\begin{align}
	\iota_{\bA} := \left( \tr{\bA}, |\dev\bA| \right),
\end{align}
where $|\cdot|$ is given by
%
$    |\bA|^2 := {1}/{2}\tr{\left( \bA^2\right)}.
$
%
For the strain-rate tensor, we denote its deviatoric part as 
\begin{align}
    \bS\bu := \dev{(\bD\bu)}
\end{align}
and note that 
%
$\tr{\bD\bu} = \nabla\cdot\bu$.
%
As a consequence of \eqref{eq:defn_eq_pressure}, $\psi_1(\bzero) = 0$ must hold. 

This characterization of isotropic rheologies allows us to narrow down the rheology inference problem to the discovery of two functions $\psi_1$ and $\psi_2$. The key physical principle of \textit{frame-indifference} then automatically holds when we represent $\Cc$ as in \eqref{eq:representation_of_C} \cite{truesdell2000}. We point out that a compressible Newtonian fluid, for which $\bsigma$ depends on $\bD\bu$ linearly, is generally written in terms of two scalar parameters, a bulk viscosity $\zeta$ and a shear viscosity $\eta$, as
\begin{align}\label{eq:compressible_newtonian_fluid}
    \bsigma = -p\bI + \zeta (\tr\bD\bu)\bI + 2 \eta\bS\bu.
\end{align}
It is straightforward to see that \eqref{eq:compressible_newtonian_fluid} is a complete characterization of all linear models with the representation \eqref{eq:representation_of_C} of $\Cc$ for which $\psi_1(\bzero) = 0$. The fact that $\psi_1 = \zeta \tr\bD\bu$ and $\psi_2 = 2 \eta$ in \eqref{eq:compressible_newtonian_fluid} establishes, for general non-Newtonian fluids,  the close relationships between $\psi_1$ and $\psi_2$ and the effective bulk and shear viscosities, respectively.   

We further restrict the functional form of the rheology with two requirements that enforce key physical and mathematical principles. Firstly, we assume the existence of a continuously differentiable and strictly convex potential of dissipation $j:\RR\times\RR_{\geq0} \to \RR$. We assume $j$, which takes the principal invariants of the strain-rate tensor $\iota_{\bD\bu}$ as input, to satisfy
\begin{equation}\label{eq:potential_of_dissipation}
\begin{split}
\frac{\partial j}{\partial (\tr{\bD\bu})}(\iota_{\bD\bu}) &= \psi_1(\iota_{\bD\bu}),\\
    \frac{\partial j}{\partial (|\bS\bu|)}(\iota_{\bD\bu}) &= \psi_2(\iota_{\bD\bu})|\bS\bu|. 
\end{split}
\end{equation}
An important consequence of this is that the Stokesian approximations to \eqref{eq:momentum_balance}, where $\bD\bu/\bD t = 0$, correspond to minima of an energy functional $\Jc:V\to\RR$ defined on a Banach space $V$, which may incorporate Dirichlet boundary conditions. 
However, this equivalence only holds when the external forces $\bff(u)$ can be represented in terms of a Fr\'echet-differentiable potential $\Phi:V\to\RR$ as 
\begin{align}
    \langle \Phi'(\bu), \bv \rangle = \int_\Omega \bff(\bu)\cdot\bv\,\dd x. 
\end{align}
Here, $\Phi'(\bu)\in V'$ denotes the Fr\'echet derivative of $\Phi$ at $\bu$ and $\langle\cdot,\cdot\rangle$ the duality pairing on $V$. Then, the energy functional $\Jc$ is given by
\begin{align}
    \Jc(\bu) = \int_\Omega \left( j(\iota_{\bD\bu}) - p \left(\nabla\cdot\bu\right) \right)\,\dd x - \Phi(\bu).
\end{align}
To show the equivalence between the minimization of $\Jc$ and the Stokesian approximation of \eqref{eq:momentum_balance} note that, using \eqref{eq:potential_of_dissipation}, the derivative of $\Jc$ in the direction of a velocity field $\bv$ is
\begin{equation}
    \begin{split}
    \langle \Jc'(\bu), \bv\rangle = \int_\Omega \left( \psi_1(\iota_{\bD\bu}) (\nabla\cdot\bv) + \psi_2(\iota_{\bD\bu}) \left( \bS\bu:\bS\bv\right) \right)\dd x  
    - \int_\Omega  p \left(\nabla\cdot\bv\right)\,\dd x - \langle\Phi'(\bu), \bv\rangle.
    \end{split}
\end{equation}
Since the necessary optimality condition $\langle \Jc'(\bu), \bv\rangle = 0$ corresponds to the weak form of \eqref{eq:momentum_balance} when $\bD\bu/\bD t = 0$, the equivalence is established. 
Moreover, if $\Phi$ is convex, then $\Jc$ is strictly convex (since we assume $j$ is strictly convex). Thus, if $\Jc$ is coercive on $V$, we can establish the existence and uniqueness of solutions to \eqref{eq:momentum_balance} whenever $\bD\bu/\bD t = 0$, see e.g.~\cite{evans2010}.

A second requirement that we enforce on the rheology is that $\Cc(\bD\bu)$ is a purely dissipative component of the Cauchy stress tensor $\bsigma$, such that
\begin{align}\label{eq:dissipative}
    \Cc(\bD\bu):\bD\bu \geq 0
\end{align}
for all velocity fields $\bu$. Inequality \eqref{eq:dissipative} is a consequence of the second law of thermodynamics \cite{morro2023}, which states that, in an isolated system, entropy must always increase. This translates into the following condition for $\psi_1$ and $\psi_2$:
\begin{align}\label{eq:dissipative_psi}
    \psi_1(\iota_{\bD\bu})(\nabla\cdot\bu) + \psi_2(\iota_{\bD\bu}) |\bS\bu|^2 \geq 0.
\end{align}
The independence of $\tr{\bD\bu}$ and $\bS\bu$ implies that \eqref{eq:dissipative_psi} holds for all velocities $\bu$ if and only if 
\begin{align}\label{eq:dissipative_psi_2}
    \psi_1(\iota_{\bD\bu})(\nabla\cdot\bu)\geq 0 \quad \text{and} \quad \psi_2(\iota_{\bD\bu}) \geq 0.
\end{align}
We may also write \eqref{eq:dissipative_psi_2} in terms of the potential of dissipation $j$ as
\begin{align}
    \frac{\partial j}{\partial (\tr{\bD\bu})}(\iota_{\bD\bu}) \tr{\bD\bu} \geq 0 \quad \text{and} \quad \frac{\partial j}{\partial( |\bS\bu|)}(\iota_{\bD\bu}) \geq 0,
\end{align}
which must hold for all velocity fields $\bu$.

\subsection{Inferring the effective shear viscosity of a fluid}\label{subsec:simplified}

We consider two simplified configurations in which the rheology inference problem is reduced to learning the dependence of $\psi_2$ on $|\bS\bu|$. That is, instead of inferring the functions $\psi_1$ and $\psi_2$, together with an equation of state for the equilibrium pressure $p$, we seek to reconstruct a function $\psi:\RR_{\geq 0}\to\RR_{\geq 0}$ given by
\begin{align}\label{eq:psi_defn}
    \psi(|\bS\bu|) := \psi_2(0, |\bS\bu|).
\end{align}
The function $\psi_2$ is proportional to the effective shear viscosity, so learning the function $\psi$ is equivalent to inferring the effective shear viscosity of a fluid. 

The first simplified configuration we consider consists of one-dimensional problems in which the fluid velocity field and the external forces can be written as
\begin{align}
    \bu(x,y,t) = \left( u(y, t), 0 \right) \quad \text{and} \quad \bff(x,y,t) = \left(f(y,t), 0\right).
\end{align}
Then, from \eqref{eq:momentum_balance}, we find that
\begin{align}
    \rho \frac{\partial u}{\partial t} - \frac{1}{2}\frac{\partial}{\partial y} \left( \psi\left(\frac{1}{2}\left|\frac{\partial u}{\partial y}\right|\right) \frac{\partial u}{\partial y} \right) = f,\label{eq:one_dim_momentum_balance}\\
    \frac{\partial}{\partial y}\left( \psi_1\left(0, \frac{1}{2}\left|\frac{\partial u}{\partial y}\right|\right) - p\right) = 0.
\end{align}
As a result, in this one-dimensional configuration, the velocity field is determined by $\psi$ (i.e.~$\psi_2$) and is independent of $\psi_1$.

The second configuration we consider is that of  incompressible fluids, where the density of the fluid $\rho$ is constant. In this case, the equation for conservation of mass simplifies to    
\begin{align}\label{eq:incompressibility}
    \nabla\cdot\bu = 0.
\end{align}
This implies that $\tr{\bD\bu} = 0$ and the rheological functions $\psi_1$ and $\psi_2$ only depend on $|\bS\bu|$. Moreover, under the assumption of incompressiblity, $\psi_1$ no longer plays a role in the dynamics of the fluid. This is a consequence of  \eqref{eq:incompressibility} being enforced with a Lagrange multiplier function $\hat{p}$. This Lagrange multiplier 
must be proportional to the hydrostatic component of $\bsigma$. Commonly, one writes
%
$    \hat{p} = - {1}/{2}\tr{\bsigma}$.
%
It then follows that $\hat{p} = p - \psi_1$. Thus, once again, the velocity field $\bu$ is independent of $\psi_1$, in the sense that knowledge of $\psi_1$ is unnecessary for computing $\bu$. For an incompressible fluid, \eqref{eq:momentum_balance} can be written as 
\begin{align}\label{eq:incompressible_momentum_balance}
    \rho \frac{\bD\bu}{\bD t} - \nabla\cdot\left(\psi(|\bS\bu|) \bS\bu\right) + \nabla \hat{p} = \bff(\bu).
\end{align}
In the two simplified setups we consider above, for the potential of dissipation $j$ given by
\begin{align}
    j(|\bS\bu|) = \int_0^{|\bS\bu|} \psi(s)s\,\dd s,
\end{align}
we can formulate \eqref{eq:one_dim_momentum_balance} and \eqref{eq:incompressible_momentum_balance} as minimization problems, granted that external forces can be written in terms of a potential and inertial terms are neglected. 
What remains to be enforced is the strict convexity of $j$. Clearly, $j$ is strictly convex if and only if the function 
%
    $s\mapsto \psi(s)s$
%
is strictly monotonically increasing for all $s\geq 0$. Finally, to ensure the dissipative nature of viscous stresses, $\psi \geq 0$ must also hold. In this way, we satisfy condition \eqref{eq:dissipative_psi_2}.

\subsection{Representing the effective shear viscosity with a neural network}\label{subsec:shear_visc_as_NN}

Under the simplified conditions from the previous section, we seek the function $\psi$ defined in \eqref{eq:psi_defn}. At this point, we introduce the additional scalar parameter $\lambda$ on which the viscosity also depends. Thus, we seek to infer the function
\begin{align}
    \psi = \psi(|\bS\bu|, \lambda).
\end{align}
In our numerical experiments, $\lambda$ represents the temperature $T$ in section \ref{sec:land_ice}, where we infer Glen's law for land ice, and the sea ice concentration $A$ in section~\ref{sec:sea_ice}, where we infer a rheology for sea ice.  

To avoid choosing an a-priori functional form for the rheology, we represent $\psi$ in terms of (feed-forward) neural networks (NNs). This allows us to exploit the excellent nonlinear function fitting capacities of NNs and make use of automatic differentiation computational modules that simplify their implementation. 
In particular, we consider two NNs
\begin{align}
	\xi:\RR\to\RR \quad \text{and} \quad \chi:\RR_{\geq0}\times \RR\to\RR_{\geq0}.
\end{align}
We use hyperbolic tangent activation functions in these NNs since they are infinitely differentiable, and are constant for large input. This latter property results in constant viscosity for strain-rates beyond the range of the data, stabilizing the model under very large strain-rates.

If we denote by $\btheta$ the vector containing the parameters (weights and biases) that characterize $\xi$ and $\chi$, we parameterize $\psi$ with the family of functions $\psi_{\btheta}$, which we define as
\begin{align}\label{eq:psi_as_NN}
	\psi_{\btheta}(|\bS\bu|, \lambda) := \exp{(\xi(\lambda))}\, \chi(\log{|\bS\bu|},\lambda).
\end{align}
In sections \ref{sec:land_ice} and \ref{sec:sea_ice} we consider fluids in which a small change in the parameter $\lambda$ leads to a large increase in the shear stress; to account for this, we apply the exponential function to $\xi$. Similarly, we aim to learn a fluid's rheology over several orders of magnitude of $|\bS\bu|$. For this reason, we precompose the first input in $\chi$ with a logarithmic function. Additionally, we wish the physical and mathematical principles stated in section \ref{subsec:isotropy} to hold for our model. For this, we require $\psi\geq 0$ and the map $s\mapsto \psi(s,\lambda)$ to be monotonically increasing. The non-negativity of $\psi$ is enforced by placing an ELU activation unit increased by one at the output $\chi$, such that $\chi\geq 0$. The monotonicity condition on $\psi(s,\lambda)s$ is weakly enforced during the optimization with a penalty term, see section \ref{subsec:optimizing_the_NN} below.

Under this configuration, the problem of inferring a rheology is reduced to that of finding the network weights $\btheta$ that minimize a given loss function. Below, in section \ref{subsec:optimizing_the_NN}, we define two loss functions that use different types of data.

\subsection{Loss functions}\label{subsec:optimizing_the_NN}

We compute the NN weights $\btheta$ that characterize the effective shear viscosity $\psi_{\btheta}$ by minimizing a loss function that quantifies the misfit between model and training data. In the simplified setting considered in this work, it suffices to generate data in steady one-dimensional scenarios. That is, we build a training dataset by finding $N$ different steady states to a one-dimensional problem. For each steady state $k = 1, 2, ..., N$, which we compute for a parameter value $\lambda^{(k)}$, we extract a stress dataset $\Sc^k$ and a velocity dataset $\Vc^k$. The stress dataset contains $M$ pairs of shear strain-rate $|\bS\bu|$ and shear stress $\sigma_{xy}$, which we denote by
\begin{align}\label{eq:defn_gamma_tau}
    \dot{\gamma} := |\bS\bu| \quad \text{and} \quad \tau := \sigma_{xy},
\end{align}
where $\sigma_{xy}$ is the off-diagonal component of $\bsigma$. In the one-dimensional setting, we have
%
$    \tau = \psi(|\dot{\gamma}|, \lambda)\dot{\gamma}
$
and $\dot{\gamma} = 1/2|\partial u/\partial y|$.
Thus, the stress dataset $\Sc^k$ and the velocity dataset $\Vc^k$ for the $k$-th steady state are given by 
\begin{align}
    \Sc^k = \{ (\dot{\gamma}^{(k)}_i, \tau^{(k)}_i)\}_{i =1}^M\:\: \text{ and } \:\:
    \Vc^k = \{ u^{(k)}_i \}_{i =1}^M.
\end{align}
The values for $\dot{\gamma}_i$, $\tau_i$, and $u_i$ are extracted by evaluating the corresponding physical quantities at points $y_i$ of the domain. The two datasets suggest two different loss functions: one stress-based and another velocity-based. The stress loss function is defined by
\begin{align}\label{eq:stress_misfit}
	\Jc_s(\btheta) := \frac{1}{N} \sum_{k = 1}^N \sum_{(\tau, \dot{\gamma})\in\Sc^{k}} \left| \log{\left(|\tau|\right)} - \log{\left(|\psi_{\btheta}(|\dot{\gamma}|, \lambda^{(k)})\,\dot{\gamma}|\right)}\right|^2.
\end{align}
The stress loss $\Jc_s$ uses a logarithmic measure of the error to account for the large differences in stress data encountered in the numerical tests. 

If we denote by $u^{(k)}_{\btheta}$ the steady-state solution $u$ to \eqref{eq:one_dim_momentum_balance} for problem $k$ in the training dataset, we can define a velocity loss function as follows:
\begin{align}
	\Jc_v(\btheta) := \frac{1}{N} \sum_{k = 1}^N \sum_{i = 1}^M \left| u_i^{(k)} - u^{(k)}_{\btheta}(y_i)\right|^2.
\end{align}
A single evaluation of $\Jc_v$ requires solving $N$ one-dimensional PDEs, i.e.~one for each of simulation in the training dataset. Computing gradients of $\Jc_v$ with respect to $\btheta$ requires differentiating though the PDE \eqref{eq:one_dim_momentum_balance}, which can be done efficiently using a combination of adjoints and AD as discussed in section~\ref{subsec:comp_setup}.

In practice, we add two penalization terms to the loss functions. Instead of minimizing the functions $\Jc_s(\btheta)$ or $\Jc_v(\btheta)$, we minimize
\begin{align}\label{eq:penalized_loss}
    \hat{\Jc}_\square(\btheta) := \Jc_\square(\btheta) + \beta_1 \|\btheta\|_{\ell_1} + \beta_2 \Pi(\btheta) \quad \text{for $\square\in\{s, v\}$.}
\end{align}
Here, $\beta_1\|\btheta\|_{\ell_1}$ penalizes the $\ell_1$ norm of the parameter vector to promote sparsification of the NN parameters, thus prioritizing simpler functions $\psi_{\btheta}$ to represent the effective shear viscosity. The function $\Pi$ is defined as
\begin{align}\label{eq:mon_penalty}
	\Pi(\btheta) := \int_{Q} \left|\min{\left\lbrace \frac{\partial}{\partial\dot{\gamma}}\left(\psi_{\btheta}(|\dot{\gamma}|, \lambda)\,\dot{\gamma}\right), 0 \right\rbrace}\right|^2\,\dd\dot{\gamma}\dd \lambda
\end{align}
and penalizes points where the derivative of the function $\dot{\gamma} \mapsto \psi(|\dot{\gamma}|, \lambda)\dot{\gamma}$ is negative over a subset $Q\subset\RR^2$ encompassing values of interest for the strain-rate $\dot\gamma$ and parameter $\lambda$. In this way, we weakly enforce the requirement that $\psi(|\dot{\gamma}|, \lambda)\dot{\gamma}$ is monotonically increasing in $\dot{\gamma}$.

\subsection{Computational aspects}\label{subsec:comp_setup}

Sections \ref{subsec:shear_visc_as_NN} and \ref{subsec:optimizing_the_NN} outline the computational framework for inferring the effective shear viscosity $\psi_{\btheta}$. Its implementation combines two Python modules, the machine learning library PyTorch \cite{paszke2019} and the finite element library Firedrake \cite{ham2023}. We use PyTorch to construct the NN-based function $\psi_{\btheta}$ and minimize the loss functions $\hat{\Jc}_s(\btheta)$ and $\hat{\Jc}_v(\btheta)$. The minimization is carried out with PyTorch's implementation of the LBFGS algorithm. 
LBFGS computes an approximation of the Hessian of the loss function with the secant method \cite{wright1999}. The gradients of $\Jc_s$ and of the penalty terms in \eqref{eq:penalized_loss} are calculated using PyTorch's automatic differentiation module. 

Gradients of $\Jc_v$ are computed with the adjoint method. In particular, the adjoint method enables efficient and accurate computations of the gradients of each of the terms inside the sum in $\Jc_v$. For computing the gradients, the boundary value problem \eqref{eq:one_dim_momentum_balance} and its adjoint equation must be solved. Thus, evaluating the gradient of $\Jc_v$ requires the solution of $2N$ PDEs, with $N$ of these being nonlinear in general. We solve these NN-based PDEs with the finite element method using Firedrake. In particular, for all one-dimensional PDEs in this paper, we use continuous piecewise-linear elements for the horizontal velocity $u$ on a mesh with $50$ cells of equal length. The required coupling between Firedrake and PyTorch is facilitated by novel functionalities that have recently been implemented in Firedrake \cite{bouziani2024}. These novel enhancements of Firedrake allow us to combine the differentiation capabilities of both modules. We also remark that, in our finite element-based approach, the PDE \eqref{eq:one_dim_momentum_balance} is a constraint to the velocity loss $\Jc_v$ whose discretization is solved accurately at each optimization step.

\section{Numerical results I: Recovering Glen's law in land ice}\label{sec:land_ice}

Ice sheets and mountain glaciers make up the large masses of ice known as land ice. These bodies of ice play an important role in Earth's climate in many ways: they regulate sea levels \cite{edwards2021}, influence ocean currents \cite{lohmann2021}, and carve out topographical features in the landscape \cite{koppes2015}, among others. Over long timescales, land ice is usually modeled as a slow-moving incompressible viscous fluid, such that the inertial terms in \eqref{eq:incompressible_momentum_balance} can be neglected. Under the condition of incompressibility, the rheology is completely defined by the effective shear viscosity function $\psi$, as explained in section \ref{subsec:simplified}.

Glen's law is the most common choice of rheological law for land ice \cite{glen1958}. It establishes a power law formula for the effective shear viscosity that depends on the temperature $T$. Using the terminology from section \ref{subsec:isotropy}, we can write Glen's law as  
\begin{align}\label{eq:glens_law}
    \psi(|\bS\bu|, T) = B(T)\left( |\bS\bu|^2 + \epsilon^2 \right)^{\frac{1-n}{2n}},
\end{align}
where $n$ is a stress exponent and $\epsilon$ a regularization parameter that sets a finite upper bound for the shear viscosity. The temperature-dependent parameter $B$ represents the fluidity of ice and is usually expressed in the form of an Arrhenius law:
\begin{align}
    B(T) := B_0 \exp{ \left(q \left(\frac{1}{T} - \frac{1}{T_{\mathrm{ref}}}  \right)\right)},
\end{align}
where $B_0$ is a pre-exponential rate factor, $q$ depends on parameters such as ice's activation energy and the universal gas constant, and $T_{\mathrm{ref}}$ is a reference temperature. Following common practice in glaciology, we set $n = 3$. For the other parameters, we set $\epsilon = 10^{-8}$, $B_0 = 1$, $q = \SI{2405}{\per\kelvin}$ and $T_{\mathrm{ref}} = \SI{273}{\kelvin}$. These quantities correspond to physically realistic values in a land ice setting under a suitable normalization.

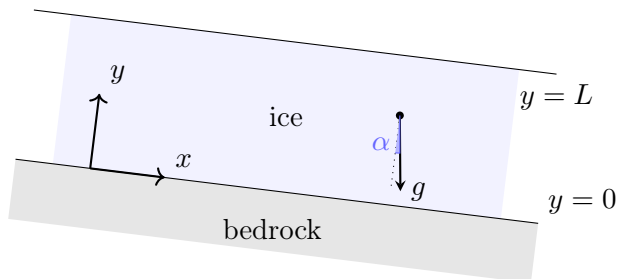
\begin{figure}
    \centering

\begin{tikzpicture}[scale=1]

\begin{scope}[rotate=-7] 
\filldraw[fill=blue!5!white,draw=none] (-0.5,2) -- (5.5,2) -- (5.5,0) -- (-0.5,0) -- (-0.5,2);
\filldraw[fill=black!10!white,draw=none] (-1,0) -- (6,0) -- (6,-.8) -- (-1,-.8) -- (-1,0);

\draw[] (-1,2) -- (6,2) node[below] {$y=L$};
\draw[] (-1,0) -- (6,0) node[above right] {$y=0$};

\node at (2.5,1) {ice};
\node at (2.5,-0.5) {bedrock};

\fill (4,1.2) circle (1.5pt);

\draw[-stealth,thick] (4,1.2) -- (4.13,0.2) node[right]{$g$};

\draw[dotted] (4,1.2) -- ++(0,-1);

\filldraw[fill=blue!30!white, draw=blue!60!white]
(4,1.2) +(0,-0.5) arc[start angle=-90,end angle=-83,radius=0.5] -- (4,1.2);
\node at (3.8,0.8) {\textcolor{blue!60!white}{$\alpha$}};

\end{scope}
\end{tikzpicture}
    \caption{Problem setup for generating a training dataset to infer Glen's law for land ice. An infinitely long slab of ice of uniform thickness $L$ slides down a bedrock inclined at an angle $\alpha$. The training dataset consists of steady solutions to this problem for different angles $\alpha$ and temperatures $T$.}
    \label{fig:pstokes_test}
\end{figure}

\subsection{Inferring Glen's law}

Here, we explore the recovery of Glen's law from shear stress and velocity measurements. To generate training data, we use a series of rheometer-like one-dimensional simulations with different parameters to expose different rheology regimes. An alternative to using several one-dimensional simulations to generate data is to use a small number of two-dimensional simulations to collect training data, as done in \cite{sunol2025}.
Our setup consists of an infinitely-long slab of ice of uniform thickness $L$ sliding down a bedrock inclined at an angle $\alpha$, see figure \ref{fig:pstokes_test}. The external forces $\bff$ correspond with gravitational forces. Since we consider all quantities but temperatures to be normalized, we set the thickness of the slab to $L = 1$ and write the gravitational acceleration vector as
%
    $\bff = (\sin{\alpha}, -\cos{\alpha}).$
%
This setup corresponds with the one-dimensional configuration examined in section \ref{subsec:simplified}. Therefore, \eqref{eq:momentum_balance} can be reduced to a single equation for the horizontal velocity $u$:
\begin{align}\label{eq:pstokes_one_dim}
    - \frac{1}{2}\frac{\partial }{\partial y}\left( \psi\left(\frac{1}{2}\left|\frac{\partial u}{\partial y}\right|, T\right)\frac{\partial u}{\partial y}\right) &= \sin{(\alpha)} \quad \text{on $(0,1)$}.
\end{align}
Equation \eqref{eq:pstokes_one_dim} is complemented with two boundary conditions: no-slip conditions at $y = 0$ and stress-free conditions at $y = 1$, such that 
\begin{align}
    u &= 0 \quad \text{at $y = 0$},\\
    \frac{\partial u}{\partial y} &= 0 \quad \text{at $y = 1$}. 
\end{align}
We generate a dataset of different velocities and stress-strain pairs by solving \eqref{eq:pstokes_one_dim} for five different slopes $\alpha = 0.01, 0.025, 0.05, 0.075, 0.1$ and three temperatures $T = 0, -10, \SI{-20}{\celsius}$ (a total of 15 configurations). For each $\alpha$ and $T$, we evaluate the solution $u$ and the shear stress $\tau$, given by \eqref{eq:defn_gamma_tau}, at ten equidistant points along $(0,1)$. We add multiplicative Gaussian noise to the stress data and additive Gaussian noise to the velocity data to mimic observational errors. In both cases, the noise has zero mean and standard deviation $\sigma_s$ and $\sigma_v$ for the stress and the velocity, respectively. For each velocity solution $u$, the additive noise for the velocity is scaled by the maximum value of the velocity $u$ in $[0,1]$, which we denote by $u_{\max}$.

\begin{table}[t]
    \centering
    \caption{Land ice problem: Average stress and velocity errors between truth and model inferred by minimizing the stress loss $\hat{\Jc}_s$ and the velocity loss $\hat{\Jc}_v$. The stress and velocity errors, denoted by $\epsilon_s$ and $\epsilon_v$, respectively, are defined in \eqref{eq:error_defn}, and are computed for error levels $\sigma_s$ and $\sigma_v$ and averaged over ten models inferred from different noise samples.}
    \label{tab:pstokes}
    \begin{tabular}{ccc|ccc}
        \multicolumn{3}{c|}{minimize $\hat{\Jc}_s$} & \multicolumn{3}{c}{minimize $\hat{\Jc}_v$} \\
        $\sigma_s$ & $\epsilon_s$ & $\epsilon_v$ & $\sigma_v/u_{\max}$ & $\epsilon_s$ & $\epsilon_v$  \\
        \midrule
        0.01 & 3.87e-05 & 6.70e-05 & 0.01 & 4.96e-05 & 2.91e-05 \\
        0.05 & 1.53e-04 & 1.06e-03 & 0.05 & 6.53e-05 & 1.01e-04 \\
        0.10 & 5.01e-04 & 9.58e-03 & 0.10 & 7.52e-04 & 2.88e-04 \\ 
    \end{tabular}
\end{table}

\begin{figure}[t]
    \centering
      \begin{subfigure}{\linewidth}
        \centering
        \includegraphics[width=0.8\linewidth]{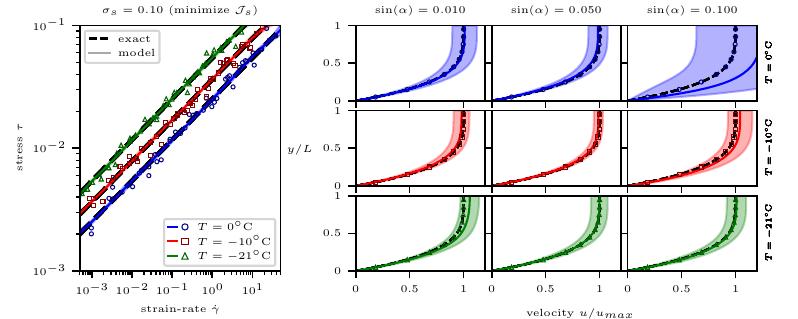}
      \end{subfigure}
      \begin{subfigure}{\linewidth}
        \centering
        \includegraphics[width=0.8\linewidth]{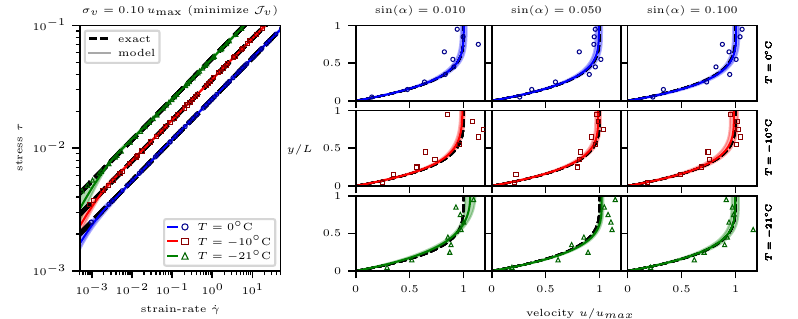}
      \end{subfigure}      
    \caption{Land ice problem: Glen's law inferred from noisy stress data (top panels) and velocity data (lower panels). For each model, we show the stress-strain relationship (left) and the velocity profiles which solve \eqref{eq:pstokes_one_dim} in the training setup (right). In each case, ten models are inferred from ten different noise samples; the results in this figure represent the mean (solid lines) and twice the standard deviation (thickness of colored regions).}
    \label{fig:pstokes}
\end{figure}

We compute the effective shear viscosity $\psi_{\btheta}$ from stress and velocity data for three different noise levels. Table \ref{tab:pstokes} presents the averaged stress and velocity errors for each noise level and each type of data. These errors, which measure the distance to the exact rheology and velocities, are defined as 
\begin{subequations}\label{eq:error_defn}
\begin{align}
    \epsilon_s &:= \frac{1}{N_{\lambda}} \sum_{i = 1}^{N_{\lambda}} \int_{\dot{\gamma}_{\min}}^{\dot{\gamma}_{\max}} \left[ \log(|\psi(|\dot{\gamma}|, T_i)\dot{\gamma}|) - \log(|\psi_{\btheta}(|\dot{\gamma}|, T_i)\dot{\gamma}|) \right]^2\,\dd \dot{\gamma}, \label{eq:stress_error_defn}\\
    \epsilon_v &:= \sum_{k = 1}^N \frac{1}{N(u^{(k)}_{\max})^2}  \int_0^L \left[ u^{(k)} - u^{(k)}(\btheta) \right]^2\,\dd y,\label{eq:vel_error_defn}
\end{align}
\end{subequations}
respectively. In \eqref{eq:stress_error_defn}, the strain-rate limits are set to $\dot{\gamma}_{\max} = 50$ and $\dot{\gamma}_{\min} = 5\times 10^{-4}$ and $N_{\lambda} = 3$ is the number of temperatures considered. The rheologies inferred from stress and velocity data with the highest error magnitudes are shown in figure \ref{fig:pstokes}. The upper set of panels correspond with models inferred from stress data with a noise level of $\sigma_s = 0.1$. Below that, we present models inferred from velocity data with noise level set to $\sigma_v/u_{\max} = 0.1$. For each of these cases, we infer ten models for different noise samples. The rheological models and velocity profiles shown in figure \ref{fig:pstokes} are mean values computed from these ten inferences. Additionally, the thickness of the colored regions around the mean values equals twice the standard deviation, visualizing the degree of uncertainty in these inferences.

 From table \ref{tab:pstokes} and figure \ref{fig:pstokes}, it is clear that Glen's law can be inferred with a high degree of accuracy from velocity data. In fact, table \ref{tab:pstokes} indicates that, for the same noise level, the stress errors of the rheological model inferred from velocity data is similar to that of the model inferred from stress data. The same, however, is not true for the velocity errors: for models inferred from stress data, the velocity errors are an order of magnitude larger than those of the models inferred from velocity data. This is clearly visible in figure \ref{fig:pstokes}, where the velocity profiles inferred from stress data (the upper set of panels) display large errors. These observations suggest that small changes in the rheology model $\psi$ lead to large changes in the resulting velocity profiles in the context of the training problem used for Glen's law. 

\subsection{Numerical minimization of the loss functions}\label{subsec:optimization_details_pstokes}

\begin{figure}[bt]
    \centering
    \includegraphics[width=0.6\linewidth]{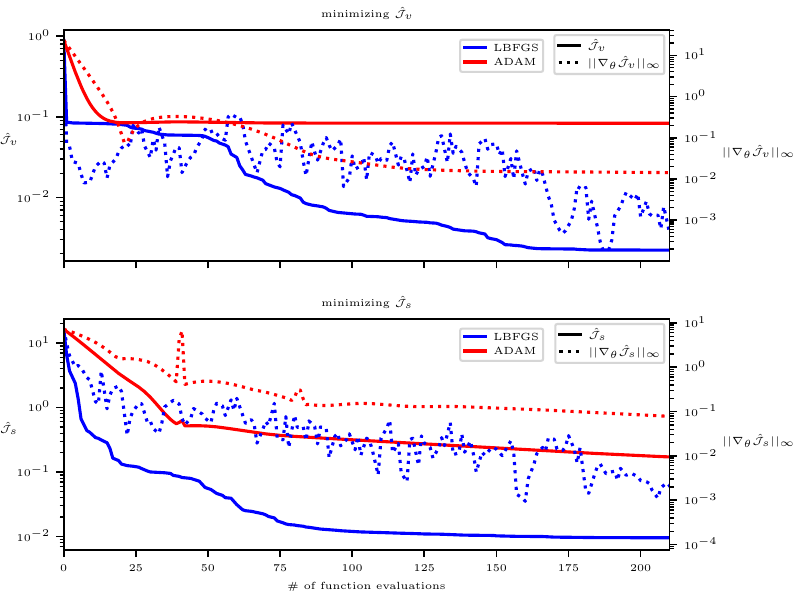}
    \caption{Convergence behavior of loss functions for Glen's law in land ice problem. The top figure shows the penalized velocity loss (solid lines) and the gradient norm (dotted lines) against the number of function evaluations when optimizing with the LBFGS (blue) and ADAM (red) algorithms. The bottom figure shows the analogue results for the penalized stress loss. Both methods use gradients computed using full batch data.}
    \label{fig:optimization_pstokes}
\end{figure}

As discussed in section \ref{subsec:comp_setup}, we minimize the loss functions $\hat{\Jc}_s$ and $\hat{\Jc}_v$ with the LBFGS algorithm implemented in PyTorch. Here, we present detailed results on the minimization of the loss functions when inferring Glen's law from stress and velocity data. The numerical behavior we observe for the problems in section~\ref{sec:sea_ice} is similar. We focus on the cases with $\sigma_s = 0.1$ and $\sigma_v/u_{\max} = 0.1$ for the stress and velocity data noise, respectively. In figure \ref{fig:optimization_pstokes}, we show the evolution of the loss functions and the $\ell^\infty$-norm of the gradients $\nabla_{\btheta}\hat{\Jc}_s$ and $\nabla_{\btheta}\hat{\Jc}_v$. To motivate the use of the LBFGS method, we also show results obtained with the ADAM algorithm, \cite{adam2014}. To compare the two algorithms fairly, we show the loss functions and their gradient norms against the number of function evaluations, instead of the number of iterations.
For ADAM, the number of function evaluations coincides with the number of iterations. However, this is not the case for LBFGS due to the linesearch method used at each iteration. For these computations, the penalty parameters are $\beta_1 = 10^{-5}$ when minimizing $\hat{\Jc}_s$ and $\beta_1 = 5\times 10^{-5}$ when minimizing $\hat{\Jc}_v$. For both loss functions, we choose $\beta_2 = 10^{14}$. For ADAM, we set the learning rate to $0.01$ and, at each iteration, we compute the gradient over the complete set of data, not over a randomized subset. We choose this learning rate value after comparing the method's performance with four other learning rates. To warm-start the minimization of the velocity loss, we perform 10 LBFGS iterations of the stress loss with noisy data. This results in an initial guess that visibly differs from the exact rheology yet is improved compared to using random NN weights.

Figure \ref{fig:optimization_pstokes} indicates that the LBFGS algorithm reduces both the velocity and stress loss functions monotonically before plateauing after about 200 function evaluations. 
A significant overall reduction in the norm of the gradient can also be observed with LBFGS despite the large variations. In contrast, ADAM (with constant step size) does not perform as well. In particular, the velocity loss $\hat{\Jc}_v$ appears to stall after a small number of iterations (see the top panel of figure~\ref{fig:optimization_pstokes}). We note that our experiments with the ADAM method have not been systematic. However, our results indicate that the LBFGS method is well-suited for the problems considered in this article.

\subsection{Two-dimensional computations}

\begin{figure}[bt]
    \centering
    \includegraphics[width=0.8\linewidth]{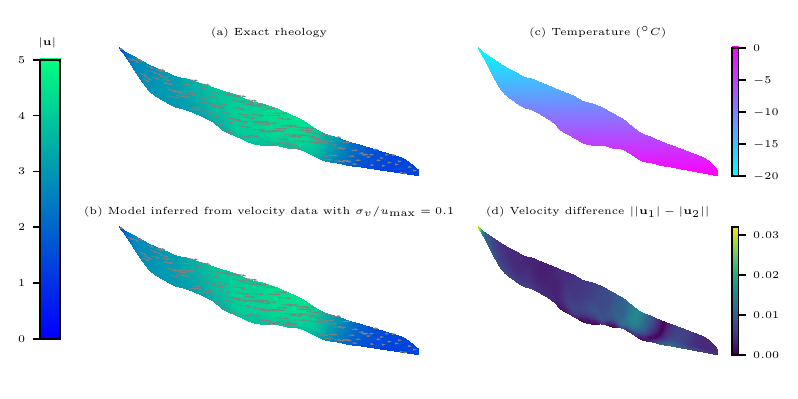}
    \caption{Land ice problem: Velocity fields computed over a longitudinal section of the Arolla glacier with (a) Glen's law \eqref{eq:glens_law} and with (b) a rheological model inferred from velocity data with a noise level of $\sigma_v = 0.1\,u_{\max}$. The color map represents the Euclidean norm of the velocity vector $|\bu|$ (non-dimensional). (c) Temperature variation, (d) difference between the two velocity fields.}
    \label{fig:two_dim_pstokes}
\end{figure}

To demonstrate the accuracy of the rheological model inferred from velocity data, we consider a more realistic two-dimensional problem. 
Extending the model to two dimensions is trivial because, as we explain in section \ref{subsec:simplified}, incompressible fluids are fully characterized by the effective shear viscosity. 
We solve the incompressible flow equations, given by \eqref{eq:incompressibility} and \eqref{eq:incompressible_momentum_balance}, over a longitudinal section of the Arolla glacier depicted in figure \ref{fig:two_dim_pstokes}; this is a common numerical test for glacier models \cite{pattyn2008}. We work with the non-dimensional values of the rheological parameters used in the dataset. 
The domain length is set to a non-dimensional value of $L = 10$ so that the resulting strain-rates are of similar orders of magnitude to those of the training dataset. 
The temperature field, which is shown in panel (c) of figure \ref{fig:two_dim_pstokes}, varies linearly with the vertical coordinate $y$, from $\SI{-20}{\celsius}$ at the highest point to $\SI{0}{\celsius}$ at the lowest point. The boundary conditions at the top boundary are stress-free conditions, such that $\bsigma\bn = 0$. At the lower boundary, we enforce no penetration conditions, $\bu\cdot\bn = 0$ and the friction law 
\begin{align}
    \bsigma_{nt} = - C |\bu_t|^{1/n-1}\bu_t,
\end{align}
where $\bsigma_{nt}$ and $\bu_t$ are the tangential components of the vectors $\bsigma\bn$ and $\bu$, respectively. The parameter $C$ is the friction coefficient, which we set to $C = 0.025$. The velocity fields computed with Glen's law and with a model inferred from velocity data with the highest error levels, $\sigma_v=0.1\,u_{\max}$, are shown in panels (a) and (b) of figure \ref{fig:two_dim_pstokes}. The two velocity fields are practically indistinguishable, with the velocity difference several orders of magnitude smaller than the velocities, as seen in panel (d). This indicates that the rheological model inferred from noisy velocity data enables accurate land ice  simulations in realistic settings.

A two-dimensional problem such as that presented in figure \ref{fig:two_dim_pstokes} contains a rich wealth of data and mimics realistic glaciological configurations better than the slab problem in figure \ref{fig:pstokes_test}. For an incompressible fluid, such a two-dimensional simulation could have been used to generate a training dataset.
This approach is followed in \cite{sunol2025}, where a shear-thinning Carreau-Yasuda model is inferred from a single two-dimensional computation in which the fluid is driven by a pressure gradient through a constriction. 
Future work could compare the effects of data errors from one and two-dimensional settings on the accuracy of the inference.

\section{Numerical results II: Inferring a rheology for sea ice}\label{sec:sea_ice}

Sea ice, which is a critical component of Earth's climate system, consists of an aggregate of ice segments, called ice floes, driven by ocean and wind drag. It is typically modeled as a continuum when studying its dynamics over spatial scales larger than $\SI{100}{km}$ \cite{feltham2008, blockley2020}. Continuum models treat sea ice like a viscous fluid that evolves over a two-dimensional horizontal plane \cite{coon1974}. Momentum balance is expressed with a thickness-averaged version of \eqref{eq:momentum_balance} that can be written as follows:
\begin{align}\label{eq:sea_ice_momentum_balance}
	\rho H \frac{D\bu}{Dt} - \nabla\cdot\bsigma = \rho_o C_o |\bu_o - \bu|(\bu_o - \bu) + \rho_a C_a |\bu_a|\bu_a.
\end{align}
In \eqref{eq:sea_ice_momentum_balance}, $H$ denotes the ice thickness. The forcing terms on the right represent ocean and atmosphere drag, which are parameterized in terms of the densities of ocean water $\rho_o$ and air $\rho_a$ and the drag coefficients $C_o$ and $C_a$. Here, we set these parameters to $\rho_o = \SI{1027}{\kg\per\cubic\meter}$, $\rho_a = \SI{1.2}{\kg\per\cubic\meter}$, $C_o = 3\times10^{-3}$ and $C_a = 10^{-3}$. Conservation of momentum must be complemented with an equation for conservation of mass, which models the evolution of the ice concentration $A$ (which takes values in $[0,1]$):
\begin{align}\label{eq:cons_mass}
	\frac{DA}{Dt} = - (\nabla\cdot\bu)A.
\end{align}

In \eqref{eq:sea_ice_momentum_balance}, the divergence of the Cauchy stress tensor $\bsigma$ represents the inner stresses arising from floe-floe interactions such as collisions, frictional contact or more complex phenomena like ridging. Sea ice is generally assumed to behave like a compressible fluid because ice floes may disperse or accumulate in different regions. As a result, for a general two-dimensional problem, inferring a rheology for sea ice requires the discovery of three different functions: the two viscosity functions $\psi_1$ and $\psi_2$ and an equation of state for the equilibrium pressure $p$. The solution to a sea ice continuum model generally consists of the velocity field $\bu$ together with two scalar fields, the ice thickness $H$ and concentration $A$.

As a first approximation, in this section we perform the rheological inference in a one-dimensional configuration in which $\psi$, the effective shear viscosity,  is the only unknown function to retrieve, as explained in section \ref{subsec:isotropy}. We consider ice floes driven by horizontal ocean and atmosphere currents, such that
\begin{align}
	\bu_o(x,y,t) = \left(  \begin{array}{c}
		u_o(y,t) \\
		0 
	\end{array} \right) \quad \text{and} \quad 	\bu_a(x,y,t) = \left(  \begin{array}{c}
		u_a(y,t) \\
		0 
	\end{array} \right).
\end{align}
We also assume that the sea ice thickness $H$ and concentration $A$ are spatially constant and, for each problem, are provided as data. Then, it is reasonable to set $\bu(x,y,t) = (u(y, t),0)$ and write \eqref{eq:sea_ice_momentum_balance} on an interval $(0,L)$ as 
\begin{align}\label{eq:sea_ice_one_dim}
    \begin{split}
        \rho H \frac{\partial u}{\partial t} - \frac{1}{2}\frac{\partial }{\partial y}\left( \psi\left(\frac{1}{2}\left|\frac{\partial u}{\partial y}\right|, A\right)\frac{\partial u}{\partial y}\right) = \rho_o C_o |u_o - u|(u_o - u) + \rho_a C_a |u_a|u_a,
    \end{split}
\end{align}
where the function $\psi:\RR_{\geq0}\times[0,1]\to\RR$ is defined as in \eqref{eq:psi_defn}. 

\begin{figure}[t]
	\centering
	\includegraphics[scale=1]{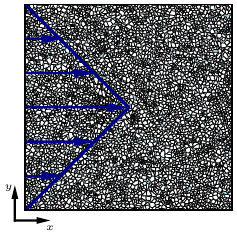}
	\caption{One-dimensional problem setup used for sea ice problems: This setup is used to generate data for both the viscous-plastic model (section~\ref{subsec:hibler}) and the DEM model (section~\ref{subsec:DEM}). In this configuration, ice floes on a periodic square patch of ocean of length $L = \SI{100}{km}$ are driven from left to right by a hat-shaped horizontal ocean current $u_o$ (blue arrows).}
	\label{fig:inference_setup}
\end{figure}

In the next two sections, we apply the approach from section \ref{sec:framework} to infer the function $\psi$ from data. In section \ref{subsec:hibler}, we recover the viscous-plastic model from data generated by solving the corresponding PDE and adding noise. Then, in section \ref{subsec:DEM}, we generate data with a discrete element method (DEM) that models polygonally-shaped ice floes. In this second scenario, the data includes significant errors due to the grid projection and averaging, and it is even uncertain whether a rheology can be found that accurately reflects the emerging flow behavior.
For both models, we generate a training dataset by finding steady states to \eqref{eq:sea_ice_one_dim} on a periodic domain of length $L = \SI{100}{\km}$. We set the atmospheric velocity $u_a$ to zero and the ocean velocity $u_o$ equal to the hat-shaped profile depicted in figure \ref{fig:inference_setup}, which is given by 
\begin{align}\label{eq:uo_triangle_profile}
	u_o(y) := U_o \left( 1 - \left| 1 - 2y/L\right| \right),
\end{align}
where $U_o$ is the maximum ocean velocity. We produce $N = 28$ steady states by solving for four different concentrations $A = 0.8, 0.85, 0.9, 0.95$ and seven ocean velocities $U_o = 0.05, 0.1, 0.25, 0.5, 1, 1.5, \SI{2}{\meter\per\second}$. In all of the computations performed in this section, the ice thickness is set to $H = \SI{2}{\meter}$. We remark that, unlike section \ref{sec:land_ice}, physical quantities are written with units in this section due to the DEM's use of dimensional quantities.

\subsection{Recovering the viscous-plastic model}\label{subsec:hibler}

The viscous-plastic model, also known as Hibler's model, is the most commonly used model for simulating sea ice dynamics \cite{hibler1979}. Practically all Earth system models use this model or its variants \cite{blockley2020}. The viscous-plastic model treats sea ice like a compressible fluid with a plastic behavior determined by an elliptical yield curve. This rheology can be fully characterized with an equation for the equilibrium pressure $p$ and the two rheological functions $\psi_1$ and $\psi_2$ introduced in section \ref{subsec:isotropy}. The two rheological functions are given by
\begin{align}
    \psi_1(\iota_{\bD\bu}, A) := \frac{p(A)}{\Delta(\iota_{\bD\bu})} \tr{\bD\bu} \quad \text{and} \quad \psi_2(\iota_{\bD\bu}, A) := \frac{p(A)}{2e^2 \Delta(\iota_{\bD\bu})},
\end{align}
where
\begin{align}
    \Delta(\iota_{\bD\bu}) := \sqrt{e^{-2}|\bS\bu|^2 + \left(\tr\bD\bu\right)^2 + \Delta_{\min}^2}.
\end{align}
Here, $\Delta_{\min}$ is a regularization parameter that leads to viscous behavior for small strain-rates, which we set to $\Delta_{\min} = \SI{2.5e-6}{\per\s}$. The parameter $e = 2$ is the eccentricity of the elliptical yield curve. 
The equilibrium pressure $p$, also known as the ice strength, is a function of the sea ice concentration $A$ and thickness $H$:
\begin{align}
    p(A) = p^\ast H \exp{(-C(1-A))}.
\end{align}
The parameters in the expression above are set to $p^\ast = \SI{2000}{\newton\per\square\meter}$ and $C = 20$. The one-dimensional model, given by \eqref{eq:sea_ice_one_dim}, is defined in terms of the effective shear viscosity $\psi$, which depends on $\psi_2$ according to \eqref{eq:psi_defn}. In the context of the viscous-plastic model, we have
\begin{align}\label{eq:hibler_rheology_one_dim}
    \psi(|\dot{\gamma}|, A) = \frac{p(A)}{2e}\frac{1}{\sqrt{|\dot{\gamma}|^2 + (e\Delta_{\min})^2}}.
\end{align}

\begin{table}[t]
    \centering
    \caption{Viscous-plastic sea ice problem: Average stress and velocity errors between truth and model inferred by minimizing the stress loss $\hat{\Jc}_s$ and the velocity loss $\hat{\Jc}_v$. The stress and velocity errors, denoted by $\epsilon_s$ and $\epsilon_v$, respectively, are defined in \eqref{eq:error_defn}, and are computed for error levels $\sigma_s$ and $\sigma_v$ and averaged over ten models inferred from different noise samples.}
    \label{tab:hibler}
    \begin{tabular}{ccc|ccc}
        \multicolumn{3}{c|}{minimize $\hat{\Jc}_s$} & \multicolumn{3}{c}{minimize $\hat{\Jc}_v$} \\
        $\sigma_s$ & $\epsilon_s$ & $\epsilon_v$ & $\sigma_v/u_{\max}$ & $\epsilon_s$ & $\epsilon_v$  \\
        \midrule
        0.01 & 2.77e-04 & 6.58e-08 & 0.01 & 4.23e-03 & 2.02e-06 \\
        0.05 & 4.40e-04 & 1.30e-07 & 0.05 & 5.55e-02 & 2.70e-05 \\
        0.10 & 6.79e-04 & 4.06e-07 & 0.10 & 4.58e-01 & 1.25e-04 \\ 
    \end{tabular}
\end{table}

\begin{figure}[t]
    \centering
      \begin{subfigure}{\linewidth}
        \centering
        \includegraphics[width=0.8\linewidth]{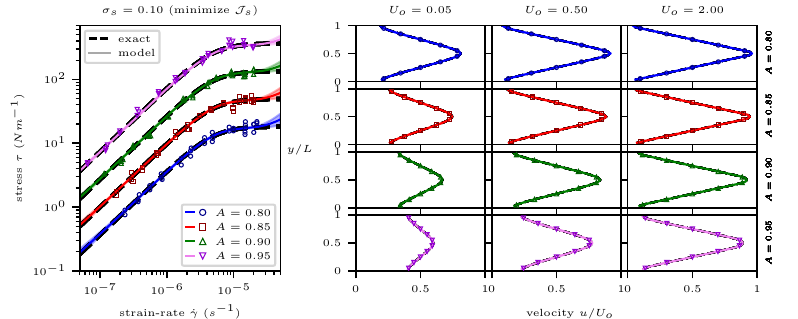}
      \end{subfigure}
      \begin{subfigure}{\linewidth}
        \centering
        \includegraphics[width=0.8\linewidth]{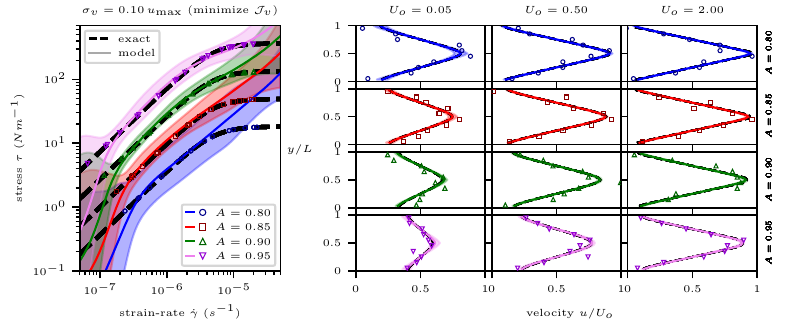}
      \end{subfigure}      
    \caption{Viscous-plastic sea ice problem: Effective shear viscosity inferred from noisy stress data (top panels) and velocity data (bottom panels). For each model, we show the stress-strain relationship (left) and the velocity profiles which solve \eqref{eq:sea_ice_one_dim} with \eqref{eq:hibler_rheology_one_dim} in the training setup (right). In each case, ten models are inferred from ten different samples of the error; the results in this figure represent the mean (solid lines) and twice the standard deviation (thickness of colored regions).}
    \label{fig:hibler}
\end{figure}

By solving \eqref{eq:sea_ice_one_dim} with $\psi$ given by \eqref{eq:hibler_rheology_one_dim}, we compute the $N = 28$ steady states and generate the training data. Following section \ref{sec:land_ice}, we add multiplicative Gaussian noise to the stress data and additive Gaussian noise to the velocity data. In both cases, the mean value is set to zero and the standard deviations to $\sigma_s$ for the stress and $\sigma_vu_{\max}$ for the velocity, where $u_{\max}$ is the maximum sea ice velocity for each steady state. 
We infer the effective shear viscosity $\psi$ from stress and velocity data for three different error magnitudes. Table \ref{tab:hibler} presents the average stress and velocity errors, which are defined in \eqref{eq:stress_error_defn} and \eqref{eq:vel_error_defn}, for three different noise levels. When computing these errors for the viscous-plastic rheology, in $\epsilon_s$ we set $\dot\gamma_{\max} = \SI{5e-8}{\per\second}$, $\dot\gamma_{\min} = \SI{5e-5}{\per\second}$, and $N_{\lambda} = 3$. We show the inferred rheologies for the highest error levels in figure \ref{fig:hibler}, where the upper and lower sets of panels correspond with models inferred from stress and velocity data, respectively. As in section \ref{sec:land_ice}, for each type of inference and noise level, we infer ten models from different noise samples. The effective shear viscosity functions shown in figure \ref{fig:hibler} thus correspond with mean values, and the thickness of the colored regions equals twice the standard deviation. 

From figure \ref{fig:hibler} and table \ref{tab:hibler}, we see that the models inferred from stress data reproduce the exact velocity profiles accurately. Conversely, the velocity data does not enable such an accurate reconstruction of the rheological model. From the lower left panel in figure \ref{fig:hibler}, it is clear that large perturbations of the rheology lead to small changes in the velocity profiles. This suggests a very different behavior to what we found in the land ice numerical test in section \ref{sec:land_ice}: small errors in the velocity data lead to large errors in the effective shear viscosity. 

\subsection{Inferring a rheology from DEM data}\label{subsec:DEM}

Discrete element methods (DEMs) resolve the dynamics of individual ice floes, including floe-scale processes such as collisions, ridging, and fracturing. Considered a promising avenue for sea ice modeling, the high computational costs involved in DEM simulations restricts their use to small scale simulations. For this reason, sea ice modules in Earth System Models rely on continuum models whose accuracy is limited by the use of phenomenological parameterizations to represent complex floe scale phenomena \cite{feltham2008}. 

To demonstrate the potential of our approach in revealing unknown rheological models from data with substantial errors, we infer an effective shear viscosity that reproduces data generated with \emph{SubZero}, a sophisticated DEM for sea ice modeling \cite{manucharyan2022}. This DEM models the dynamics of irregularly-shaped polygonal ice floes that interact through collisions, see figure \ref{fig:inference_setup}. More complex processes, such as fracturing and ridging, are modeled in this DEM. However, we deactivate these features when generating the data used in the inference. The numerical results presented here are similar to those in \cite{dediego2025}, and we refer to that work for further information on the rheology inferred from DEM data and several tests proving the rheology's ability to generalize beyond the training dataset. The results presented here differ from \cite{dediego2025} in that two rheological models are inferred from stress and velocity data independently (in \cite{dediego2025}, the model inferred from velocity data uses the model inferred from stress data as an initial guess in the optimization). Moreover, unlike \cite{dediego2025}, the NN $\chi$ in 
\eqref{eq:psi_as_NN} is precomposed with a logarithm. We remark that these changes do not result in substantially different models.  

\begin{figure}[t]
    \centering
    \includegraphics[width=0.8\linewidth]{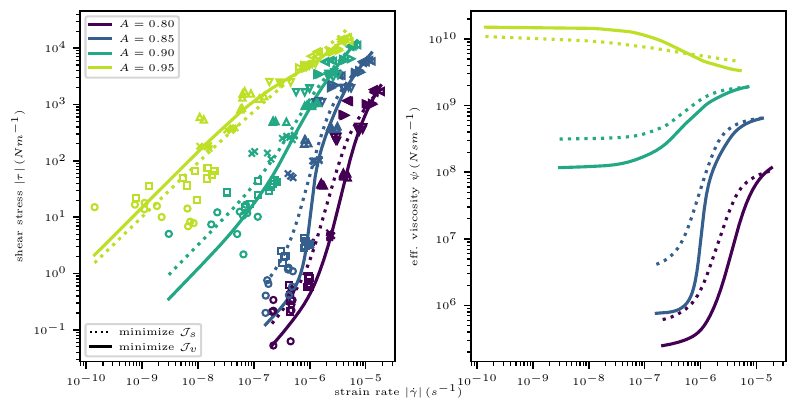}
    \caption{DEM sea ice problem: Shown on the left is the shear strain-rate/stress relationship for different concentrations (solid) learned from the DEM data (markers). Shown on the right is the effective shear viscosity $\psi_{\btheta}$ for different concentrations. We show the rheological models computed by minimizing the stress loss $\hat{\Jc}_s$ (dotted lines) and the velocity loss $\hat{\Jc}_v$ (solid lines).}
    \label{fig:subzero_rheology}
\end{figure}

\begin{figure}[t]
    \centering
    \includegraphics[width=0.95\linewidth]{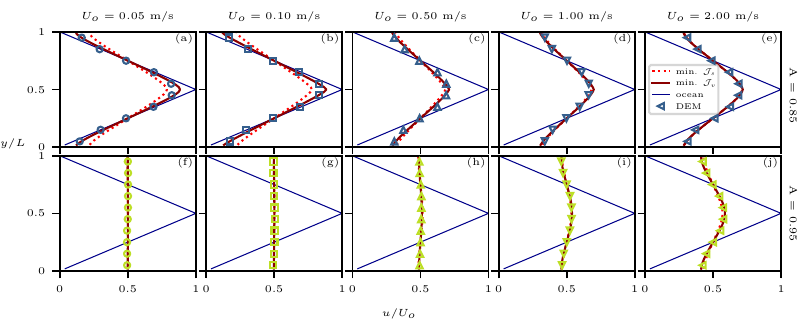}
    \caption{DEM sea ice model: Velocity profiles computed with the rheological model inferred from minimizing the velocity loss $\hat{\Jc}_v$ (solid red lines) and stress loss $\hat{\Jc}_s$ (red dotted lines) over the training dataset. We also show the DEM data (markers) and ocean velocity profiles (blue lines). Each row corresponds with a sea ice concentration $A$ and each column with a maximum ocean velocities $U_o$.}
    \label{fig:subzero_vel}
\end{figure}

Two rheological models inferred from DEM data are presented in figure \ref{fig:subzero_rheology}, one computed by minimizing the stress loss $\hat{\Jc}_s$ and the other the velocity loss $\hat{\Jc}_v$. The left panel displays the two stress-strain relationships, together with the DEM data points. The right panel presents the effective shear viscosity $\psi$ of the two models. Velocity profiles computed by solving the training dataset problems with the two models are shown in figure \ref{fig:subzero_vel}. We also include the values of the (unpenalized) loss functions $\Jc_s$ and $\Jc_v$ in table \ref{tab:subzero}. The rheological models predicted with both stress and velocity data reveal a transition from a shear-thickening to a shear-thinning viscosity as the sea ice concentration $A$ increases. The model inferred from stress data achieves a high degree of accuracy in reproducing the velocity profiles of the training dataset (see the red dotted line in figure \ref{fig:subzero_vel}). However, some discrepancies between model and velocity data arise for low ocean velocities $U_o$ and $A = 0.85$; this is unsurprising, given the substantial errors present in the DEM stress data due to interpolation from Lagrangian to Eulerian data, and the necessary spatial and temporal averaging. The model inferred from velocity data, however, corrects this discrepancy, yielding an impeccable match between model prediction and velocity training data. Further tests presented in \cite{dediego2025} demonstrate the models' accuracy in reproducing DEM data outside of the training dataset.

\begin{table}[t]
    \centering
    \caption{DEM sea ice model: Values of the (unpenalized) stress and velocity loss functions $\Jc_s$ and $\Jc_v$ for the two rheological models that result from minimizing $\hat{\Jc}_s$ and $\hat{\Jc}_v$. }
    \label{tab:subzero}
    \begin{tabular}{c|cc}
        & minimize $\hat{\Jc}_s$ & minimize $\hat{\Jc}_v$ \\
        \midrule
        $\Jc_s$ & $0.991$ & $1.93$ \\
        $\Jc_v$ & $2.24\times 10^{-4}$ & $6.88\times 10^{-5}$
    \end{tabular}
\end{table}

\section{Discussion and future work}

\subsection{Extensions to multiple external parameters}
In this work, we only consider a single external parameter (temperature for land ice and concentration for sea ice). However, the rheologies of complex systems often depend on additional parameters; for example, the size and thickness of ice floes are known to influence the overall sea ice rheology \cite{wilchinsky2006, feltham2008, herman2022}. Including multiple external parameters increases the size of the training dataset and thus the computational cost of fitting the NN parameters.
Since in this work (and most related work cited in the introduction), we assume that our training data results from (possibly expensive) simulations rather than physical observations or measurements, we are able to, at least in principle, generate training data at will. In the presence of several external parameters, one could generate (and use during training) only a small number of selected combinations of external parameters, instead of a tensor grid of parameters as in this work.   Which parameter combinations to generate data for could be chosen through random sampling, or systematically such that the simulations expose the rheological behavior in different regimes. Selective data generation that optimally informs the rheology would be an interesting topic of future research. It could be based on sensitivities of the loss function, or other criteria to identify optimal training data motivated from recent advances in active learning or experimental design.

\subsection{Sensitivity of inferred models to data errors}
The numerical results presented in sections \ref{sec:land_ice} and \ref{sec:sea_ice} indicate very different sensitivities of the models to errors in the data. For the land ice model, in section \ref{sec:land_ice}, we find that adding small noise to the stress data results in large deviations from the true velocity profile, as observed in the top set of panels in figure \ref{fig:pstokes}. This should be contrasted with the lower set of panels in figure \ref{fig:pstokes}, where noise in the velocity data translates into small errors in the effective shear viscosity. Interestingly, for the sea ice model in section \ref{subsec:hibler}, we observe the opposite effect. In figure \ref{fig:hibler} we see that noise in the stress data results in small errors in the velocity computed with the model; however, noise in the velocity data implies large errors in the effective shear viscosity. 

These results suggest that the
effects of noise on the accuracy of a rheology inferred from velocity data are very much problem-dependent. In the land and sea ice problems explored in sections \ref{sec:land_ice} and \ref{subsec:hibler}, respectively, future research should investigate the causes of such different effects of noise; advances in this direction may help establish guidelines to infer unknown rheologies (i.e.~what data should be used) and estimate the errors in these inferences. It is unclear to what extent these differences are caused by the underlying rheological models (Glen's law and the viscous-plastic model) or by the presence of ocean drag in the sea ice model. We remark that ocean drag is an additional nonlinearity in the governing PDE that has large effects on the velocity profiles. Under steady conditions and no atmospheric forcing, conditions under which we train the sea ice model, inner stresses balance ocean drag. Since ocean drag depends on the difference between sea ice and ocean velocity quadratically, small changes in the rheology imply even smaller changes in the sea ice velocity. In fact, under certain circumstances, we can expect velocity solutions of the viscous-plastic model to be close to insensitive to changes in the rheology. For example, for fast ocean velocities or small sea ice concentrations, the sea ice velocity approaches a free-drift regime where $u\approx u_o$, see the velocity profiles in figure \ref{fig:hibler} for $U_o = \SI{2}{\meter\per\second}$. Under free-drift conditions, weakening the material (i.e.~reducing the effective viscosity) translates into $u$ further approaching $u_o$, which will not be noticeable if $u\approx u_o$. 

\subsection{Inference of the convex potential of dissipation}
The framework presented in section \ref{sec:framework} suggests an alternative approach to rheology inference: Instead of learning the effective shear viscosity $\psi$, learn the potential of dissipation $j$. The potential $j$ is introduced in section \ref{subsec:isotropy} to establish conditions under which the governing equations are well-posed. To this end, we enforce the convexity of $j$ or, equivalently, the monotonicity of $s\mapsto\psi(s, \lambda)s$ when the fluid is one-dimensional or incompressible. Parameterizing the convex potential $j$ with input convex NNs, which have demonstrated accurate interpolation performance for convex models \cite{amos2017}, would remove the need for the monotonicity penalty. Moreover, for two-dimensional compressible fluids, this approach  ensures the existence of a convex potential of dissipation (in two dimensions, unlike the one-dimensional or incompressible fluid case, an underlying potential of dissipation does not exist in general). Input convex NNs are used in \cite{parolini2025} to represent convex effective viscosities successfully. Despite many rheological models being convex, this is not always the case; see e.g.~the effective shear viscosities of polymer solutions in \cite{laun1984}.  

\section{Conclusions}

This article presents a framework for inferring rheological models from data. By characterizing isotropic rheologies governed by an underlying convex potential of dissipation, we ensure isotropic frame-indifference and the well-posedness of the inferred models. We apply this framework to infer the effective shear viscosity function for one-dimensional configurations and incompressible fluids. By parameterizing the effective shear viscosity with a neural network, we set up a computational framework capable of training the network with both stress and velocity data. We minimize the misfit between velocity data and the model using a PDE-constrained optimization approach that solves the governing equations with the finite element method.

We apply our method to infer the effective shear viscosities of land and sea ice models. To understand the accuracy of our method, we first infer two known rheologies, Glen's law for land ice and the effective shear viscosity of the viscous-plastic model for sea ice. Our numerical results demonstrate our method's robustness in recovering the rheological model under large errors. Interestingly, we find large differences between the land and sea ice problems in the effects of noise on the accuracy of the recovered models. Finally, we infer a rheology that reproduces the velocity fields that result from Lagrangian ice floe simulations computed with a discrete element method. No rheology is known to reproduce the dynamics of this synthetic sea ice system. However, our method is able to construct a nonlinear effective shear viscosity that transitions from shear-thickening to a shear-thinning behavior as the sea ice concentration increases. This model reproduces the averaged Lagrangian ice floe velocity fields accurately in one-dimensional settings.

\section*{Acknowledgement}
The authors were supported by the Multidisciplinary University Research Initiatives (MURI) Program, Office of Naval Research (ONR) grant \#N00014-19-1-242. G.S.~also appreciates support from the US National Science Foundation (NSF) under \#2343866, \#2411229 and \#2411349.

\section*{Data availability}

The data and scripts required to reproduce the results in this article are contained in the following repository: 
\url{https://bitbucket.org/gonzalogddiego/rheologyinference}.

\bibliography{bibliography}

@article{jenkins1983,
  title={A theory for the rapid flow of identical, smooth, nearly elastic, spherical particles},
  author={Jenkins, James T and Savage, Stuart B},
  journal={Journal of Fluid Mechanics},
  volume={130},
  pages={187--202},
  year={1983},
  publisher={Cambridge University Press}
}

@book{ranalli1995,
  title={Rheology of the {E}arth},
  author={Ranalli, Giorgio},
  year={1995},
  publisher={Springer Science \& Business Media}
}

@book{larson2013,
  title={Constitutive equations for polymer melts and solutions},
  author={Larson, Ronald G},
  year={2013},
  publisher={Butterworth-Heinemann}
}

@article{hu2024,
author = {Jiashun Hu  and Johann Rudi  and Michael Gurnis  and Georg Stadler },
title = {Constraining {E}arth’s nonlinear mantle viscosity using plate-boundary resolving global inversions},
journal = {Proceedings of the National Academy of Sciences},
volume = {121},
number = {28},
pages = {e2318706121},
year = {2024},
}

@article{akerson2025,
  title={{Learning constitutive relations from experiments: 1. PDE-constrained optimization}},
  author={Akerson, Andrew and Rajan, Aakila and Bhattacharya, Kaushik},
  journal={Journal of the Mechanics and Physics of Solids},
  pages={106128},
  year={2025},
  publisher={Elsevier}
}

@article{kontogiannis2025,
  title={Learning rheological parameters of non-{N}ewtonian fluids from velocimetry data},
  author={Kontogiannis, Alexandros and Hodgkinson, Richard and Reynolds, Steven and Manchester, Emily L},
  journal={Journal of Fluid Mechanics},
  volume={1011},
  pages={R3},
  year={2025},
  publisher={Cambridge University Press}
}

@article{saadat2022,
  title={Data-driven selection of constitutive models via rheology-informed neural networks {(RhINNs)}},
  author={Saadat, Milad and Mahmoudabadbozchelou, Mohammadamin and Jamali, Safa},
  journal={Rheologica Acta},
  volume={61},
  number={10},
  pages={721--732},
  year={2022},
  publisher={Springer}
}

@article{mahmoud2024,
  title={Unbiased construction of constitutive relations for soft materials from experiments via rheology-informed neural networks},
  author={Mahmoudabadbozchelou, Mohammadamin and Kamani, Krutarth M and Rogers, Simon A and Jamali, Safa},
  journal={Proceedings of the National Academy of Sciences},
  volume={121},
  number={2},
  pages={e2313658121},
  year={2024},
  publisher={National Academy of Sciences}
}

@article{thakur2024,
title = {{ViscoelasticNet: {A} physics informed neural network framework for stress discovery and model selection}},
journal = {Journal of Non-{N}ewtonian Fluid Mechanics},
volume = {330},
pages = {105265},
year = {2024},
author = {Thakur, Sukirt and Raissi, Maziar and Ardekani, Arezoo M.},
}

@article{wang2025,
  title={Deep learning the flow law of {A}ntarctic ice shelves},
  author={Wang, Yongji and Lai, Ching-Yao and Prior, David J and Cowen-Breen, Charlie},
  journal={Science},
  volume={387},
  number={6739},
  pages={1219--1224},
  year={2025},
  publisher={American Association for the Advancement of Science}
}

@article{lardy2025,
title = {Inferring viscoplastic models from velocity fields: {A} physics-informed neural network approach},
journal = {Journal of Non-{N}ewtonian Fluid Mechanics},
volume = {346},
pages = {105512},
year = {2025},
author = {Martin Lardy and Sham Tlili and Simon Gsell},
}

@article{reyes2021,
  title={Learning unknown physics of non-{N}ewtonian fluids},
  author={Reyes, Brandon and Howard, Amanda A and Perdikaris, Paris and Tartakovsky, Alexandre M},
  journal={Physical Review Fluids},
  volume={6},
  number={7},
  pages={073301},
  year={2021},
}

@article{bruer2022,
  title={Inferring ice sheet damage models from limited observations using {CRIKit}: the {C}onstitutive {R}elation {I}nference {T}oolkit},
  author={Bruer, Grant and Isaac, Tobin},
  journal={arXiv preprint arXiv:2204.09748},
  year={2022}
}

@article{lennon2023,
  title={Scientific machine learning for modeling and simulating complex fluids},
  author={Lennon, Kyle R and McKinley, Gareth H and Swan, James W},
  journal={Proceedings of the National Academy of Sciences},
  volume={120},
  number={27},
  pages={e2304669120},
  year={2023},
  publisher={National Academy of Sciences}
}

@article{parolini2025,
  title={Structure-preserving neural networks in data-driven rheological models},
  author={Parolini, Nicola and Poiatti, Andrea and Ven{\'e}, Julian and Verani, Marco},
  journal={SIAM Journal on Scientific Computing},
  volume={47},
  number={1},
  pages={C182--C206},
  year={2025},
  publisher={SIAM}
}

@article{dediego2025,
  title={Non-{N}ewtonian viscous fluid models with learned rheology accurately reproduce {L}agrangian sea ice simulations},
  author={de Diego, Gonzalo G and Stadler, Georg},
  journal={Physical Review Fluids},
  volume={11},
  number={1},
  pages={013301},
  year={2026},
  publisher={APS}
}

@article{sunol2025,
      title={Learning constitutive models and rheology from partial flow measurements}, 
      author={Alp M. Sunol and James V. Roggeveen and Mohammed G. Alhashim and Henry S. Bae and Michael P. Brenner},
      year={2025},
      journal={preprint, arXiv:2510.24673},
      archivePrefix={arXiv},
      url={https://arxiv.org/abs/2510.24673}, 
}

@article{coon1974,
  title={Modeling the pack-ice as an elastic-plastic material},
  author={Coon, M.D. and Maykut, G.S. and Pritchard, R.S. and Rorthrock, D.A. and Thorndike, A.S.},
  journal={{AIDJEX} Bulletin},
  volume={24},
  pages={1-105},
  year={1974}
}

@article{hibler1979,
  title={A dynamic thermodynamic sea ice model},
  author={Hibler III, WD},
  journal={Journal of physical oceanography},
  volume={9},
  number={4},
  pages={815-846},
  year={1979}
}

@article{wilchinsky2006,
  title={Modelling the rheology of sea ice as a collection of diamond-shaped floes},
  author={Wilchinsky, Alexander V and Feltham, Daniel L},
  journal={Journal of non-Newtonian fluid mechanics},
  volume={138},
  number={1},
  pages={22--32},
  year={2006},
  publisher={Elsevier}
}

@article{feltham2008,
  title={Sea ice rheology},
  author={Feltham, Daniel L},
  journal={Annu. Rev. Fluid Mech.},
  volume={40},
  pages={91-112},
  year={2008},
}

@article {blockley2020,
  title={The future of sea ice modeling: {W}here do we go from here?},
  author={Blockley, Ed and Vancoppenolle, Martin and Hunke, Elizabeth and Bitz, Cecilia and Feltham, Daniel and Lemieux, Jean-Fran{\c{c}}ois and Losch, Martin and Maisonnave, Eric and Notz, Dirk and Rampal, Pierre and others},
  journal={Bulletin of the American Meteorological Society},
  volume={101},
  number={8},
  pages={E1304--E1311},
  year={2020}
}

@article{herman2022,
  title={Granular effects in sea ice rheology in the marginal ice zone},
  author={Herman, A.},
  journal={Philosophical Transactions of the Royal Society A},
  volume={380},
  number={2235},
  pages={20210260},
  year={2022},
  publisher={The Royal Society}
}

@article{toppaladoddi2025,
  title={A viscous continuum theory of sea ice motion based on stochastic floe dynamics},
  author={Toppaladoddi, S},
  journal={Journal of Fluid Mechanics},
  volume={1014},
  pages={A6},
  year={2025},
  publisher={Cambridge University Press}
}

@article{manucharyan2022,
  title={{SubZero}: A sea ice model with an explicit representation of the floe life cycle},
  author={Manucharyan, Georgy E and Montemuro, Brandon P},
  journal={Journal of Advances in Modeling Earth Systems},
  volume={14},
  number={12},
  pages={e2022MS003247},
  year={2022},
  publisher={Wiley Online Library}
}

@article{jop2006,
  title={A constitutive law for dense granular flows},
  author={Jop, Pierre and Forterre, Yo{\"e}l and Pouliquen, Olivier},
  journal={Nature},
  volume={441},
  number={7094},
  pages={727--730},
  year={2006},
  publisher={Nature Publishing Group UK London}
}

@book{gurtin1982,
  title={An introduction to continuum mechanics},
  author={Gurtin, Morton E},
  year={1982},
  publisher={Academic press}
}

@book{truesdell2000,
  title={An introduction to the mechanics of fluids},
  author={Truesdell, Clifford and Rajagopal, Kumbakonam Ramamani},
  year={2000},
  publisher={Springer Science \& Business Media}
}

@article{rajagopal2006,
  title={On implicit constitutive theories for fluids},
  author={Rajagopal, Kumbakonam Ramamani},
  journal={Journal of Fluid Mechanics},
  volume={550},
  pages={243--249},
  year={2006},
  publisher={Cambridge University Press}
}

@article{kamrin2012,
  title={Nonlocal constitutive relation for steady granular flow},
  author={Kamrin, Ken and Koval, Georg},
  journal={Physical review letters},
  volume={108},
  number={17},
  pages={178301},
  year={2012},
  publisher={APS}
}

@book{morro2023,
  title={Mathematical modelling of continuum physics},
  author={Morro, Angelo and Giorgi, Claudio},
  volume={314},
  year={2023},
  publisher={Springer}
}

@article{wright1999,
  title={Numerical Optimization},
  author={Wright, Stephen and Nocedal, Jorge},
  journal={Springer Science},
  volume={35},
  number={67-68},
  pages={7},
  year={1999}
}

@article{adam2014,
  title={Adam: A method for stochastic optimization},
  author={Kingma Diederik P. and Ba, Jimmy},
  journal={arXiv preprint arXiv:1412.6980},
  volume={1412},
  number={6},
  year={2014}
}

@inproceedings{amos2017,
  title={Input convex neural networks},
  author={Amos, Brandon and Xu, Lei and Kolter, J Zico},
  booktitle={International conference on machine learning},
  pages={146--155},
  year={2017},
  organization={PMLR}
}

@manual{ham2023,
  title        = {Firedrake User Manual},
 author       = {David A. Ham and Paul H. J. Kelly and Lawrence Mitchell and Colin J. Cotter and Robert C. Kirby and Koki Sagiyama and Nacime Bouziani and Sophia Vorderwuelbecke and Thomas J. Gregory and Jack Betteridge and Daniel R. Shapero and Reuben W. Nixon-Hill and Connor J. Ward and Patrick E. Farrell and Pablo D. Brubeck and India Marsden and Thomas H. Gibson and Miklós Homolya and Tianjiao Sun and Andrew T. T. McRae and Fabio Luporini and Alastair Gregory and Michael Lange and Simon W. Funke and Florian Rathgeber and Gheorghe-Teodor Bercea and Graham R. Markall},
  organization = {Imperial College London and University of Oxford and Baylor University and University of Washington},
  year         = {2023},
  doi          = {10.25561/104839},
}

@article{bouziani2024,
  title={Differentiable programming across the {PDE} and Machine Learning barrier},
  author={Bouziani, Nacime and Ham, David A and Farsi, Ado},
  journal={arXiv preprint arXiv:2409.06085},
  year={2024}
}

@article{alhashim2025,
  title={Control of flow behavior in complex fluids using automatic differentiation},
  author={Alhashim, Mohammed G and Hausknecht, Kaylie and Brenner, Michael P},
  journal={Proceedings of the National Academy of Sciences},
  volume={122},
  number={8},
  pages={e2403644122},
  year={2025},
  publisher={National Academy of Sciences}
}

@book{evans2010,
  title={Partial differential equations},
  author={Evans, Lawrence C},
  volume={19},
  year={2022},
  publisher={American Mathematical Society}
}

@article{paszke2019,
  title={Pytorch: An imperative style, high-performance deep learning library},
  author={Paszke, Adam and Gross, Sam and Massa, Francisco and Lerer, Adam and Bradbury, James and Chanan, Gregory and Killeen, Trevor and Lin, Zeming and Gimelshein, Natalia and Antiga, Luca and others},
  journal={Advances in neural information processing systems},
  volume={32},
  year={2019}
}

@incollection{glen1958,
	author = {Glen, J.W.},
	title = {The flow law of ice: a discussion of the assumptions made in glacier theory, their experimental foundation and consequences},
	booktitle = {Physics of the Movement of Ice: Symposium at Chamonix 1958},
	publisher = {Int. Assoc. Hydrol. Sci.},
	address = {Wallingford, UK},
	year = {1958},
	pages = {171-183}
}

@Article{pattyn2008,
AUTHOR = {Pattyn, F. and Perichon, L. and Aschwanden, A. and Breuer, B. and de Smedt, B. and Gagliardini, O. and Gudmundsson, G. H. and Hindmarsh, R. C. A. and Hubbard, A. and Johnson, J. V. and Kleiner, T. and Konovalov, Y. and Martin, C. and Payne, A. J. and Pollard, D. and Price, S. and R\"uckamp, M. and Saito, F. and Sou\v{c}ek, O. and Sugiyama, S. and Zwinger, T.},
TITLE = {Benchmark experiments for higher-order and full-Stokes ice sheet models ({ISMIP–HOM})},
JOURNAL = {The Cryosphere},
VOLUME = {2},
YEAR = {2008},
NUMBER = {2},
PAGES = {95--108},
}

@book{greve2009,
  title={Dynamics of ice sheets and glaciers},
  author={Greve, R. and Blatter, H.},
  year={2009},
  publisher={Springer Science \& Business Media}
}

@article{edwards2021,
  title={Projected land ice contributions to twenty-first-century sea level rise},
  author={Edwards, Tamsin L and Nowicki, Sophie and Marzeion, Ben and Hock, Regine and Goelzer, Heiko and Seroussi, H{\'e}l{\`e}ne and Jourdain, Nicolas C and Slater, Donald A and Turner, Fiona E and Smith, Christopher J and others},
  journal={Nature},
  volume={593},
  number={7857},
  pages={74--82},
  year={2021},
  publisher={Nature Publishing Group UK London}
}

@article{koppes2015,
  title={Observed latitudinal variations in erosion as a function of glacier dynamics},
  author={Koppes, Mich{\'e}le and Hallet, Bernard and Rignot, Eric and Mouginot, J{\'e}r{\'e}mie and Wellner, Julia Smith and Boldt, Katherine},
  journal={Nature},
  volume={526},
  number={7571},
  pages={100--103},
  year={2015},
  publisher={Nature Publishing Group UK London}
}

@article{lohmann2021,
  title={Risk of tipping the overturning circulation due to increasing rates of ice melt},
  author={Lohmann, Johannes and Ditlevsen, Peter D},
  journal={Proceedings of the National Academy of Sciences},
  volume={118},
  number={9},
  pages={e2017989118},
  year={2021},
  publisher={National Academy of Sciences}
}

@article{Guazzelli18,
  title={Rheology of dense granular suspensions},
  author={Guazzelli, {\'E}lisabeth and Pouliquen, Olivier},
  journal={Journal of Fluid Mechanics},
  volume={852},
  pages={P1},
  year={2018},
  publisher={Cambridge University Press}
}

@article{laun1984,
  title={Rheological properties of aqueous polymer dispersions},
  author={Laun, Hans Martin},
  journal={Die Angewandte Makromolekulare Chemie: Applied Macromolecular Chemistry and Physics},
  volume={123},
  number={1},
  pages={335--359},
  year={1984},
  publisher={Wiley Online Library}
}
\bibliographystyle{abbrv}

\end{document}